\documentclass[letterpaper, 10 pt, conference]{ieeeconf}  % Comment this line out
\IEEEoverridecommandlockouts    
\usepackage{comment}
\usepackage{cite}
\newcommand{\eat}[1]{}
\usepackage{booktabs}
\usepackage{color}
\usepackage [autostyle, english = american]{csquotes}
\MakeOuterQuote{"}
\usepackage{xcolor}
\let\labelindent\relax %Otherwise throw an error since IEEE also defines indent
\usepackage{enumitem}

\usepackage{multicol}
\usepackage{multirow}

\usepackage{mathtools}
    
\usepackage{amsmath}
\usepackage{amstext}% assumes amsmath package installed
\usepackage{amssymb}
\usepackage{amsfonts}
\usepackage{float}

\usepackage{amsthm}  %for makinf assumptions bold
\usepackage[normalem]{ulem} % for strikethrough

\newtheorem{remark}{Remark}

\newtheorem{problem}{Problem}
\usepackage{subfig}
\usepackage{caption}
\usepackage{todonotes}

\makeatletter

\newcommand{\Rmnum}[1]{\expandafter\@slowromancap\romannumeral #1@}
\usepackage{savesym}

\usepackage{algorithm}
\usepackage{algorithmicx} % <===========================================
\usepackage{algpseudocode} %
\savesymbol{AND}
\usepackage[group-separator={,},group-minimum-digits={3}]{siunitx}

\usepackage{graphicx} % for pdf, bitmapped graphics files
\usepackage{epsfig} % for postscript graphics files

\usepackage{times} % assumes new font selection scheme installed
\usepackage{amsmath} % assumes amsmath package installed
\usepackage{amssymb}  % assumes amsmath package installed
\usepackage{comment}
\makeatletter
\let\NAT@parse\undefined
\makeatother
\usepackage{hyperref}
\hypersetup{
   colorlinks=true,
    linkcolor= blue,
    allcolors=blue,
    citecolor = blue,
    filecolor=black,      
    urlcolor=blue,
    }
\usepackage{mathrsfs}

\newcommand{\tf}{t_f^{*}}

\usepackage[symbol]{footmisc}

\title{\LARGE \bf
Minimally Disruptive Cooperative Lane-change Maneuvers}
\author{ Behdad Chalaki$^{1*}$, Vaishnav Tadiparthi$^{1*}$, Hossein Nourkhiz Mahjoub$^1$, Jovin D'sa$^1$, Ehsan Moradi-Pari$^1$,\\ 
Andres S. Chavez Armijos$^2$, Anni Li$^2$, and Christos G. Cassandras$^2$
\thanks{$^*$Both authors contributed equally. $^1$B. Chalaki, V. Tadiparthi, H. Nourkhiz Mahjoub, E. Moradi-Pari, and Jovin D'sa are with Honda Research Institute-US (HRI-US) Ann Arbor, MI
48103 USA (email:{behdad chalaki; vaishnav\_tadiparthi; hossein\_nourkhizmahjoub; emoradipari; jovin\_dsa}@honda-ri.com) $^2$A. S. Chavez Armijos, A. Li, and C. G. Cassandras are with the Division of Systems Engineering and Center for Information and Systems
Engineering, Boston University, Brookline, MA 02446 (email:{aschavez;
anlianni; cgc}@bu.edu).}
}

\begin{document}

\maketitle
\thispagestyle{empty}
\pagestyle{empty}

\begin{abstract} 
A lane-change maneuver on a congested highway could be severely disruptive or even infeasible without the cooperation of neighboring cars.
However, cooperation with other vehicles does not guarantee that the performed maneuver will not have a negative impact on traffic flow unless it is explicitly considered in the cooperative controller design. 
In this letter, we present a socially compliant framework for cooperative lane-change maneuvers for an arbitrary number of CAVs on highways that aims to interrupt traffic flow as minimally as possible. 
Moreover, we explicitly impose feasibility constraints in the optimization formulation by using reachability set theory, leading to a unified design that removes the need for an iterative procedure used in prior work.
We quantitatively evaluate the effectiveness of our framework and compare it against previously offered approaches in terms of maneuver time and incurred throughput disruption.
\end{abstract}

\section{Introduction}
%\todoall{related works from previous phases, motivation, some new papers.} 
%\todoall{@Hossein Cleaning up the introduction with the references. 
%@Vaishnav:  \sout{Feasibility, condensing previous approaches}.
%@Behdad: \sout{making a flowchart for our approaches (experiment to condense other approaches);simulation setup;} 
%@together: Simulation results }

%Literature review 

% Behdad: removed  AVChallenges_2
%\PARstart{D}{espite} great advances in the field of autonomous driving during almost the last two decades, this technology is yet to be realized on a large-scale due to some of its challenging aspects, such as strict safety requirements \cite{hussain2018autonomous}. 
\PARstart{D}{espite} great advances in the field of autonomous driving over the last two decades, strict safety requirements have prevented it to be realized on a large-scale \cite{hussain2018autonomous}.
Guaranteeing safety could be very difficult, especially in high-speed scenarios, or if there are multiple interacting agents in the scene. Introducing the concept of connectivity is a potential strong remedy here, due to its capability for enhancing the situational awareness of autonomous vehicles (AVs) by providing them with external streams of information \cite{CAVs_1, CAVs_2}. In addition, this technology makes AVs capable of being cooperative with each other by sharing their real-time state information. Therefore, Connected AVs (CAVs) have a notably higher capability, compared to AVs, to handle safety-critical scenarios. In addition to safety, this capability helps to improve driving efficiency and comfort as well  \cite{CAVs_3}.
%Behdad: removed  CAVs_5
%Different aspects of designing a system of multiple CAVs has been studied in the literature \cite{CAVs_4, CAVs_6}. 

Among different scenarios, highway driving, due to its relatively high speed, is one of the most attractive yet challenging ones in which CAVs have a great potential to help in increasing safety and efficiency. Different studies exist in the literature that try to tackle a highway autonomous driving problem, such as an autonomous car following design \cite{ACF_1, ACF_2}. Automating a lane change maneuver, which is an intuitively more challenging task due to its higher dimensionality compared to a longitudinal maneuver design, has also gained attention from the research community.
Many research attempts on autonomous lane change have been published in the literature, either as an advisory system \cite{ulbrich2015towards, tariq2022slas}
%schmidt2021probabilistic} 
to check the feasibility of a lane change, or as a motion planner \cite{nilsson2016lane,luo2016dynamic,he2021rule}. Although the autonomous lane-change maneuver exhibits promising results, without the cooperation of other vehicles the maneuver may become infeasible. 
To address the infeasibility issue, various research investigated cooperation among several CAVs to perform a lane change maneuver \cite{li2017optimal,li2018balancing,katriniok2020nonconvex}. Yet, most studies ignored the negative effects of the lane-change maneuver on the surrounding vehicles. Specifically, Wang et al. \cite{wang2022ego} showed that restricting an AV's negative impacts on neighboring vehicles, can results in improved traffic performance.
%Several studies considered cooperation among several CAVs to perform a lane change maneuver \cite{li2017optimal,li2018balancing,katriniok2020nonconvex} in order to address the infeasibility issue. 
%Most of these studies, however, completely ignored the negative impact of the lane-change maneuver on the surrounding vehicles. 

In a series of our previous studies, we have focused on designing a cooperative framework for a group of CAVs to perform a lane-change maneuver in a highway scenario by considering and reducing the negative effects of this maneuver on the cooperative agents \cite{chen2020cooperative,chavez2022sequential,armijos2022cooperative}.
We address the scenario depicted in Fig. \ref{fig:Cooperative lane change scneario}, in which a CAV must pass a slow uncontrolled vehicle and make a lane change to the fast lane by cooperating with a pair of CAVs on the fast lane. This task was divided into two phases: positioning the cooperating vehicles in an appropriate and safe formation only by altering their longitudinal locations, and then performing a lane-change maneuver to change the ego vehicle's lateral position.

%Contrbution of this paper and comparison with other related work
In this letter, we extend our previous research in multiple directions. First, our framework leverages cooperating with an arbitrary number of CAVs on the fast lane by considering the global impact of the maneuver. Secondly, we utilize the developments in reachability set theory to explicitly impose feasibility constraints in the optimization problem itself.
This leads to our unified formulation for coordination which removes the need for an iterative procedure used in \cite{chavez2022sequential} and thus improves upon our previous methods. To the best of our knowledge, our design advances the state of the art by allowing the ego vehicle to cooperate with more than one pair of vehicles to further reduce highway traffic disruption.
% on a highway.  

%Organization of the paper

The remaining sections of this letter are structured as follows. In Section \ref{sec:Prelim}, we briefly discuss the modeling framework's prerequisites. In Section \ref{sec:Model}, we present a summary of earlier work, followed by an elaborate description of our proposed formulations. In Section \ref{sec:simulation}, we conduct simulations to empirically demonstrate the efficacy of our approach. 
The final section \ref{sec:Conclusion} contains concluding remarks.
% In the final section, \ref{sec:Conclusion}, we provide concluding remarks.

% \section{Modeling Framework}
%\todoall{COPY FROM DOC SENT TO BU, ALSO REUSE STUFF FROM PREVIOUS PAPERS.}
%\todoall{Define dynamics and set of cooperating CAVs}
%\todoall{Define feasible sets for final position and speed}
%\todoall{Maybe explaining the steps for computing the final time based on the previous BU paper}
%\todoall{List all the assumptions such as no communication delay, there is a low-level tracker which can track the trajectories, we do not consider other non-CAV participants, this is ongoing effort}

\section{Preliminaries} \label{sec:Prelim}
\begin{figure}
    \centering
\includegraphics[width=0.75\linewidth]{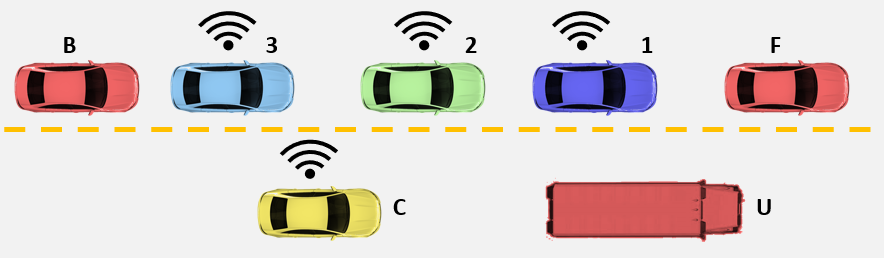}  
\caption{A scenario of three potential CAVs on the fast lane, and a CAV $C$ which is behind a slow uncontrolled vehicle $U$ and aims at performing a lane change.}
%\caption{A scenario of three potential CAVs (vehicles 1,2, and 3) on the fast lane (left lane), and a CAV (vehicle $C$) which is behind a slow uncontrolled vehicle (vehicle $U$) and aims at performing a lane-change.}
    \label{fig:Cooperative lane change scneario}
\end{figure}

We consider motion planning of CAVs on a highway to cooperatively enable an ego CAV to pass a slow uncontrolled vehicle (see Fig. \ref{fig:Cooperative lane change scneario} for a CAV group of size four). 
This lane change maneuver can be decomposed into longitudinal and lateral phases. First, cooperating CAVs and the ego CAV adjust their longitudinal position such that there is a safe and feasible gap for the ego CAV to perform a lane change. Then, the ego CAV performs an optimal lane-change maneuver. In this letter, we focus solely on the first component, i.e., the longitudinal motion planning of CAVs.

 % \subsection{Preliminaries}
 Let $C$ be the ego CAV which aims at performing a lane change to the fast lane, and passes a slow uncontrolled vehicle denoted by $U$. The set of all cooperative CAVs on the fast lane at time $t\in\mathbb{R}$, which are in the communication range of CAV $C$, is given by $\mathcal{S}(t)$. The set of cooperative CAVs on the fast lane can be located between potential human-driven vehicles $F$ and $B$ from front and back, respectively. For simplicity in notation, let $\mathcal{S}(t)=\{1,\dots, m\}$, where $1$ is the CAV farthest ahead and $m$ is the last.% CAV in the group.

We model the longitudinal dynamics of CAV $i\in \mathcal{S}(t)$ as a double integrator, i.e.,
\begin{align}
\begin{aligned}\label{eq:dynamics}
\dot{x}_i(t)=v_i(t),\\
\dot{v}_i(t)=u_i(t),
\end{aligned}\end{align}
where $x_{i}(t)$, $v_{i}(t)$, and
$u_{i}(t)$ denote position, speed, and acceleration at $t\in\mathbb{R}$, respectively, and $\mathbf{x}_{i}(t)=\left[p_{i}(t), v_{i}(t)\right]^\top$ be the state of the CAV $i$ at time $t$.  Let $t_0$ denote the time at which the longitudinal component of the maneuver is triggered, while $t_f$ is the time when this component is completed. Without loss of generality thereafter we consider $t_0=0$ to simplify the notation.

For each CAV $i\in\mathcal{S}(t)$ the control input and speed are bounded by 
\begin{equation}\label{uconstraint}
    u_{\min}\leq u_i(t)\leq u_{\max},
\end{equation}
\begin{equation}\label{vconstraint}
    0\leq v_{\min}\leq v_i(t)\leq v_{\max},
\end{equation}
where $u_{\min},u_{\max}$ are the minimum and maximum control inputs and $v_{\min},v_{\max}$ are the minimum and maximum speed limit, respectively. 

To guarantee rear-end safety between CAV $i\in \mathcal{S}(t)$ and a preceding CAV $k\in\mathcal{S}(t)$, we impose the following speed-dependent constraint,
\begin{gather}
 p_k(t)-p_i(t)\geq \epsilon + \varphi~ v_i(t), %\underbrace{\epsilon + \varphi~ v_i(t)}_{\delta_i(t)},
\end{gather}
%where $\delta_i(t)$ is the safe distance, while
where $\epsilon$ and $\varphi\in\mathbb{R}_{>0}$ are the standstill distance and reaction time respectively. 

In contrast to approaches that neglect the negative effects of a lane change on the surrounding traffic \cite{he2021rule,li2018balancing, katriniok2020nonconvex,luo2016dynamic,nilsson2016lane,li2017optimal}, we borrow a metric called disruption from \cite{armijos2022cooperative} to explicitly consider the adverse impact that the maneuver might have on the fast-lane traffic. For any vehicle $i\in{\mathcal{S}}(t)$, the total disruption at time $t>0$ is denoted by $D_i(t)$, and given by
\begin{subequations}  
\begin{align}
    &D_i(t) = \gamma_x \Delta^x_i(t) + \gamma_v \Delta^v_i(t),\label{eq:delta}\\ 
    &\Delta^x_i(t) = (x_i(t) - (x_i(0)+v_i(0)\cdot t))^2,\label{eq:deltaX}\\
    &\Delta^v_i(t) = (v_i(t) - v_d)^2,\label{eq:deltaV}
\end{align}
\end{subequations}
where $\Delta^x_i(t)$ and $\Delta^v_i(t)$ are the position and flow disruptions. We define the disruption metric using not only the expected change in the final positions of the CAVs in the fast lane but also the change in their final speed. The position disruption metric $\Delta^x_i(t)$ is a measure of the discrepancy in the terminal positions of the cooperating vehicle in the fast lane when compared to its position had it not cooperated at all and cruised with its initial speed over $[0,t]$. The flow disruption $\Delta^v_i(t)$ is defined as the speed deviation of the vehicle at time $t$ from a desired flow speed, $v_d$.

To make the dimensions consistent in \eqref{eq:delta}, the weight factors $\gamma_x$ and $\gamma_v$ are defined as follows 
\begin{subequations}  
\begin{align}
    &\gamma_x = \dfrac{\gamma}{(\max(v_{\max}-v_0, v_{\min}-v_0)\cdot t_{avg})^2},\\
    &\gamma_v = (1-\gamma)\cdot\max(v_{\max}-v_{d}, v_{\min}-v_{d})^2,
\end{align}
\end{subequations}
where $\gamma\in[0,1]$ is a tuning parameter to place more emphasis on position or flow disruption, while $t_{avg}$ is the desired average time to complete the longitudinal component of the maneuver. In the next section, we briefly describe the general approaches in the previous work \cite{chavez2022sequential, armijos2022cooperative} and elaborate on the proposed formulations.

\section{Modeling Framework}\label{sec:Model}

\subsection{Previous Formulations}
We will compare the performance of our proposed approaches primarily against two formulations developed in the recent past (\cite{chavez2022sequential, armijos2022cooperative}) to tackle the problem of a cooperative lane-change maneuver on a highway. 
%{\color{red} For convenience of notation, let us refer to the techniques presented in \cite{chavez2022sequential} and \cite{armijos2022cooperative} as {ITS-CONF} and {ITS-JNL} respectively.}

In both \cite{chavez2022sequential} and  \cite{armijos2022cooperative}, a multi-step iterative approach is presented:
\begin{enumerate}[leftmargin=*, label={Step \arabic*} ]%[label=Step \arabic{enumi}:,ref=Step \arabic{enumi}, wide=0 pt]%[leftmargin=*, label={Step \arabic*} ]
    \item \label{step1} CAV $C$ determines an optimal terminal maneuver time $t_f^{*}$ and control input (acceleration/deceleration) $u_C^* (t)$, $t \in [0, t_f]$ that minimizes a given objective function, while satisfying vehicle dynamics and safety constraints with the slow vehicle in front of it (i.e., vehicle $U$ in Fig. \ref{fig:Cooperative lane change scneario}).    
    \item \label{step2} An optimal pair of CAVs $(i^*, i^* + 1)$ is identified from
    % the set of cooperative CAVs 
    $\mathcal{S}(\tf)$ in the fast lane.
    \item \label{step3} In case no feasible pair is identified, the terminal time $\tf$ is relaxed and the series of optimization problems is re-solved from \ref{step1}.
    \item \label{step4} A planning algorithm may be executed to determine the optimal trajectory that minimizes an energy cost subject to dynamics and terminal state constraints. 
\end{enumerate}
 The major differences between the two approaches in \cite{armijos2022cooperative} and \cite{chavez2022sequential} are listed below. 

\begin{figure}[!h]
    \centering
    \includegraphics[width = \linewidth]{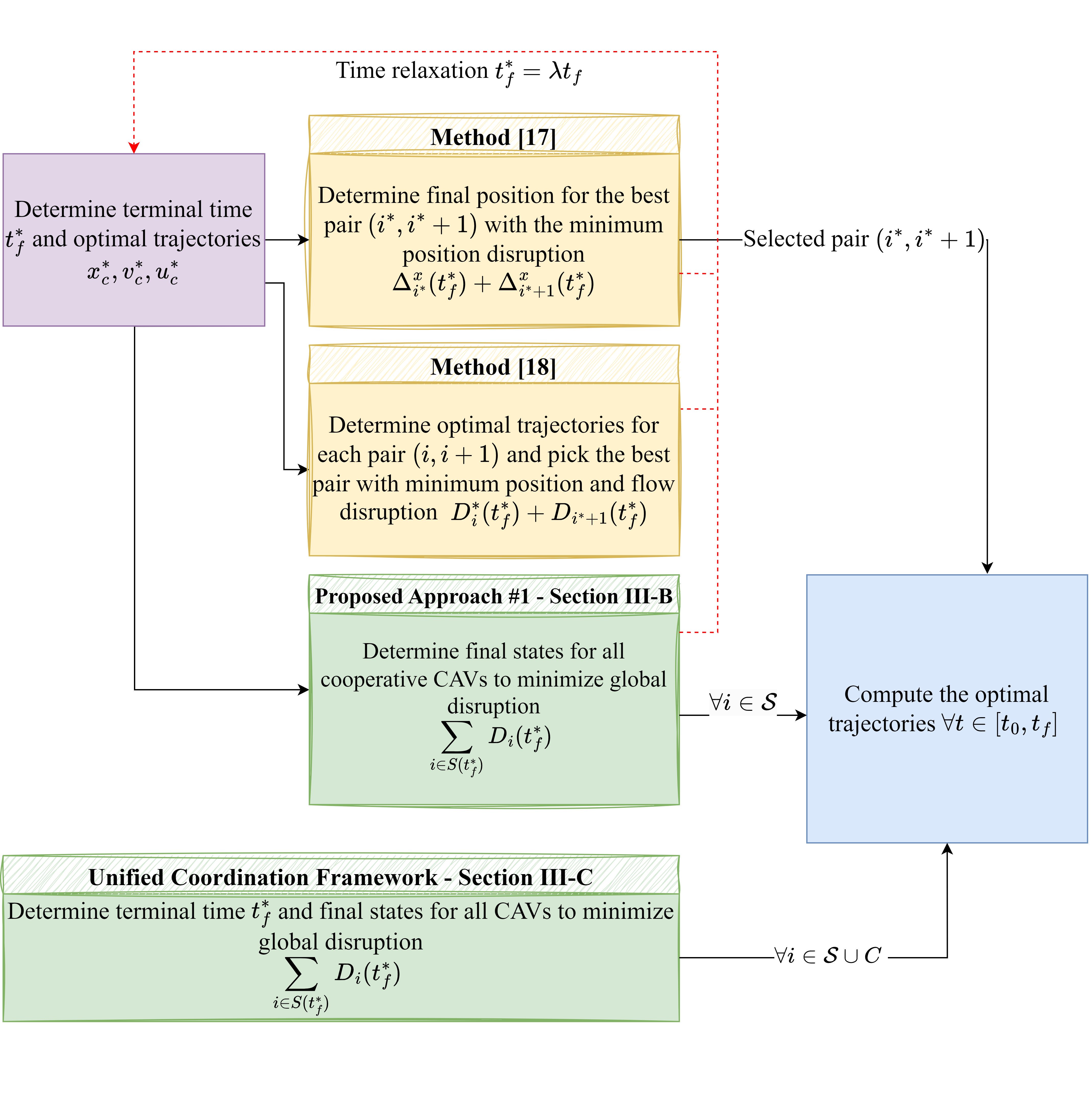}
    \caption{Block diagrams of different approaches for longitudinal positioning of CAVs in a cooperative lane change.}
    \label{fig:blocks}
\end{figure}

In \ref{step1} (purple block in Fig. \ref{fig:blocks}), for both approaches, the objective function consists of a weighted sum of the time required for the maneuver and the energy cost for CAV $C$ over the same time. 
In \cite{chavez2022sequential}, the terminal speed $v_C(\tf)$ is constrained to be close to some given $v_d$ (i.e., desired speed of fast lane) while in \cite{armijos2022cooperative}, this is as a terminal cost rather than a hard constraint as in \cite{chavez2022sequential}.
%This optimization problem determines the terminal position $x_C (\tf)$ as well.

In \ref{step2} (yellow blocks in Fig.\ref{fig:blocks}),  the optimal pair is chosen such that it minimizes a measure of disruption incurred by the cooperating vehicles. In \cite{chavez2022sequential}, the disruption metric is limited to position (i.e., $\Delta^x_{i}(\tf) +\Delta^x_{i+1} (\tf)$ from \eqref{eq:deltaX}) and the decision variables also include the terminal positions of the cooperative pair. 
However, in \cite{armijos2022cooperative}, the disruption metric additionally accounts for deviation in speed (i.e., $D_{i}(\tf) +D_{i+1} (\tf)$ from \eqref{eq:delta}). 

In \ref{step3} (red dotted line in Fig. \ref{fig:blocks}), the minimum disruption found has to be below a specified threshold $D_{th}$; otherwise,
the terminal time is relaxed as $t'_f = \lambda \tf$ (where $\lambda >1$) and another iteration takes place from \ref{step1}.

\ref{step4} (blue block in Fig. \ref{fig:blocks}) is not required for approach \cite{armijos2022cooperative} as it is combined with \ref{step2}, in which 
the planned trajectory is found in conjunction with the optimal pair of CAVs that minimizes the combined disruption. However, it is a necessary component for the approach in \cite{chavez2022sequential}.

The key observations to note from these two methods are:
\begin{enumerate}[leftmargin=*]
    \item The first optimization problem solved to find the terminal time $t_f^{*}$ in \ref{step1} does not account for any CAVs in $\mathcal{S}(.)$. %the vehicles in the fast lane.

    \item An iterative process is followed to determine the feasibility of the obtained terminal time. 
    Depending on the initial conditions of the traffic scenario and the desired speeds, this approach may be severely time-consuming due to the need for multiple iterations.%, as we will demonstrate empirically. 

    \item The disruption measured only accounts for the pair of cooperating vehicles in the fast lane ($i, i+1$), but not the vehicles behind the selected pair, i.e., CAVs $j > i+ 1$, $j \in \mathcal{S}(\cdot)$.
    %A holistic view of cooperation required to conduct the maneuver would account for any disruption incurred by these vehicles too. 
\end{enumerate}
In the next sections, we introduce our framework to address the aforementioned limitations. Fig. \ref{fig:blocks} summarizes all these different approaches.  

% \todoall{Would it be better to include the flowcharts from these papers?}
% \todobehdad{If we have space, it might be good to also include a simplified flow-chart (not from the paper, cause those were very confusing) describing this procedure in a high-level}

\subsection{Cooperating with more CAVs on the fast lane}\label{PA1}
The first change that we desire to introduce is the principle of collaboration with \textit{all} CAVs in the cooperative set of CAV $C$.
% First, we seek to revise the formulation to consider collaboration with more than a pair of CAVs in the fast lane. 
We shall illustrate this approach with the formulation developed for an arbitrary set of cooperating CAVs denoted by $\mathcal{S}(t_f^*)=\{1,\dots,m\}$, and CAV $C$ aims at performing a lane change given a final time $t_f^*$ and its final states (see Fig.\ref{fig:Cooperative lane change scneario}). 
%, i.e., step 1 is intact.
Considering \ref{step1} is intact, the solution of this section determines the terminal states of the cooperative CAVs on the fast lane and a location into which CAV $C$ merges (top green block shown in Fig. \ref{fig:blocks}). We define global disruption as the cost incurred by all the cooperating CAVs in the swarm block due to the lane-change maneuver initiated by CAV $C$. Note that, the optimal pair of CAVs which CAV $C$ merges between is a decision variable that needs to be computed along with other decision variables such that they minimize the \textbf{global} disruption of CAVs in the fast lane. 

\begin{comment}
\begin{remark}
The optimal pair of CAVs which CAV $C$ merges between is a decision variable that needs to be computed along with other decision variables such that they minimize the \textbf{global} disruption of CAVs in the fast lane. 
\end{remark}
\end{comment}

\subsubsection{Safety constraints among CAVs on the fast lane}
Since we aim to find terminal states of CAVs on the fast lane, we impose rear-end safety constraints for CAVs in the fast lane only at time $t_f^*$, and let a trajectory planner (\ref{sec:trajPlanner}) %(optimal control problem) 
ensure safety during $t\in[0, t_f^*)$. Thus, we have the following rear-end safety constraints. 
\begin{subequations}\label{eq:ConstraintsFastLane}
\begin{align}
   &x_F(t_f^*) - x_1(t_f^*)\geq \varphi v_1(t_f^*)+\epsilon, \label{eq:F,1}\\
   &\forall j\in\mathcal{S}(t)\setminus\{m\}:\nonumber\\ 
   &x_j(t_f^*)-x_{j+1}(t_f^*)\geq\varphi v_{j+1}(t_f^*)+\epsilon, \\
   &x_m(t_f^*) - x_B(t_f^*)\geq \varphi {v_B(t_f^*)}+\epsilon \label{eq:constWithB}.
   \end{align}
\end{subequations}
Note that the rationale for considering \eqref{eq:constWithB} is that CAV $m$ should not take an action that disturbs vehicle $B$ (potentially a HDV), so it considers the safety constraint with the vehicle behind it as well. 

%\begin{remark}
%We consider that states of vehicles $B$ and $F$ at time $t_f^*$ i.e., $\mathbf{x}_B(t_f^*)=[x_B(t_f^*), v_B(t_f^*)]$, $\mathbf{x}_F(t_f^*)=[x_F(t_f^*), v_F(t_f^*)]$, come from a human behavior prediction module.
%\end{remark}
\subsubsection{Feasible set of final states}
As our optimization formulations seek to derive the terminal states for all the vehicles in the set of cooperative vehicles, it is crucial to determine their feasible set.
Given the double integrator dynamics in \eqref{eq:dynamics} for all CAVs $i \in \mathcal{S}(\tf)$, we apply the geometry analysis tools presented in \cite{haddad2021curious, kurzhanski1997ellipsoidal} to determine a parametric representation for this convex set, also called the \textit{reachable} set.
In the presence of control input constraints \eqref{eq:dynamics}, define the following two parameters:
\begin{subequations}
    \begin{align}
        \mu =\frac{u_{\max}-u_{\min}}{2}, \quad 
        \nu = \frac{u_{\max}+u_{\min}}{2}.
    \end{align}
\end{subequations}
The implicit form of the bounding curves of the reachable set at time $t$ for a double integrator are parabolic in state and given by:
% Using analysis described in \cite{kurzhanski1997ellipsoidal, haddad2021curious}, the reachable set of a double integrator system \eqref{eq:dynamics} for CAV $i\in {S}(t)$ at time $t$ with control input constraint \eqref{uconstraint} can be analytically computed as a convex set 
\begin{subequations} 
    \allowdisplaybreaks
\begin{align}
    p^{\rm upper}(x_i, v_i,t) =&-\frac{t^2}{2} + \frac{1}{4} \left(\frac{v_i - v_i(0) - \nu t}{\mu} + t \right)^2 \notag\\    
    &- \frac{x_i - x_i(0) - t v_i(0) - \nu \frac{t^2}{2}}{\mu} ,\\
    p^{\rm lower}(x_i, v_i,t) =&~\frac{t^2}{2} -\frac{1}{4} \left(\frac{-v_i + v_i(0) + \nu t}{\mu} + t \right)^2 \notag\\
    &- \frac{x_i - x_i(0) - t v_i(0) - \nu \frac{t^2}{2}}{\mu}.
\end{align}
\end{subequations}

% where 

% \todoall{DO WE NEED EXAMPLE OF FEASIBILITY SET? VT: I THINK NOT.}

% For example, the reachable set for a double integrator with initial conditions $x_0 = 0$ and $v_0 = 0$ at time $t = 0.5$ is shown in Fig. \ref{fig:reach_set}.
% \begin{figure}[!h]
%     \centering
%     \includegraphics[width = 0.8\linewidth]{Figures/reach.eps}
%     \caption{Illustration of Feasibility Set}
%     \label{fig:reach_set}
% \end{figure}

The reachable set can now be established as
\begin{equation}\label{eq:reachable set}
    p^{\rm upper}(x_i, v_i,t) \leq 0 \;  \cup \;
     p^{\rm lower}(x_i, v_i,t) \geq 0.
\end{equation}
This constraint can be expressed in the decision variables as required and added to the list of constraints in the final optimization problem for terminal states of CAV $i$. 
%\todoall{\sout{VAISHNAV: Maybe adding more materials(text/figure) on this section. }}

\subsubsection{Safety constraint of CAV $C$ with corresponding vehicles on the fast lane} %based on the target location}
Focusing on cases that CAV $C$ does not merge in front of vehicle $F$ or behind vehicle $B$, which are uncontrolled vehicles, we need to consider safety constraints based on the target location of CAV $C$ and cooperative CAVs on the fast lane. In what follows, we first provide the safety constraints for a sample case shown in Fig. \ref{fig:Cooperative lane change scneario}. Starting with listing safety constraints based on the location of CAV $C$ in \eqref{eq:CAV1ALL}-\eqref{eq:CAV3ALL}, in \eqref{eq:ALLOR} we show that several of these constraints are disjunctive meaning they cannot be satisfied at the same time. Next, we formalize these constraints by incorporating binary variables in \eqref{eq:CAV1BINARY}, and extend it to the arbitrary set of cooperating CAVs in \eqref{eq:ConstraintsFastLaneCAVC}.
%%%%%%

Case 1: CAV $C$ merges between CAV $1$ and vehicle $F$,
\begin{subequations}\label{eq:CAV1ALL}
\begin{align}
   &x_F(t_f^*) - x_C(t_f^*)\geq \varphi v_C(t_f^*)+\epsilon,\label{eq:fFrontc}\\
   &x_C(t_f^*) - x_1(t_f^*)\geq \varphi v_1(t_f^*)+\epsilon. \label{eq:CFront1}
\end{align}
\end{subequations}

Case 2: CAV $C$ merges between CAV $2$ and CAV $1$,
\begin{subequations}
\begin{align}
   &x_1(t_f^*) - x_C(t_f^*)\geq \varphi v_C(t_f^*)+\epsilon,\label{eq:1FrontC}\\
   &x_C(t_f^*) - x_2(t_f^*)\geq \varphi v_2(t_f^*)+\epsilon.\label{eq:CFront2}
\end{align}
\end{subequations}

Case 3: CAV $C$ merges between CAV $3$ and CAV $2$,
\begin{subequations}
\begin{align}
   &x_2(t_f^*) - x_C(t_f^*)\geq \varphi v_C(t_f^*)+\epsilon,\label{eq:2FrontC}\\
   &x_C(t_f^*) - x_3(t_f^*)\geq \varphi v_3(t_f^*)+\epsilon.\label{eq:CFront3}
\end{align}
\end{subequations}

Case 4: CAV $C$ merges between vehicle $B$ and CAV $3$,
\begin{subequations}  \label{eq:CAV3ALL}
\begin{align}
   &x_3(t_f^*) - x_C(t_f^*)\geq \varphi v_C(t_f^*)+\epsilon,\label{eq:3FrontC}\\
   &x_c(t_f^*) - {x_B(t_f^*)}\geq \varphi {v_B(t_f^*)}+\epsilon.\label{eq:CFrontb}
\end{align}
\end{subequations}
We need to always satisfy \eqref{eq:fFrontc} and \eqref{eq:CFrontb}; however, only one of the constraints in each of the pairs \eqref{eq:CFront1} and \eqref{eq:1FrontC}, pair  \eqref{eq:CFront2} and \eqref{eq:2FrontC}, and pair  \eqref{eq:CFront3} and \eqref{eq:3FrontC} needs to be satisfied, which can be written  as 
\begin{subequations} \label{eq:ALLOR}
\begin{align}
\eqref{eq:CFront1}~~\texttt{OR}~~\eqref{eq:1FrontC},\label{eq:1ORC}\\
\eqref{eq:CFront2}~~\texttt{OR}~~\eqref{eq:2FrontC},\label{eq:2ORC}\\
\eqref{eq:CFront3}~~\texttt{OR}~~\eqref{eq:3FrontC}\label{eq:3ORC}.
\end{align}
\end{subequations}
Constraints \eqref{eq:1ORC}, \eqref{eq:2ORC}, and \eqref{eq:3ORC} are disjunctive constraints due to the $\texttt{OR}$ statement.
Moreover, they determine if CAV $C$ merges in front or behind a cooperative CAV in the fast lane.
We convert each disjunctive constraint into two separate constraints by introducing a binary variable $B_i\in\{0,1\}$ for all $i\in \mathcal{S}(\tf)$ and a sufficiently large number $M\in\mathbb{R}_{\geq 0}$ \cite{grossmann2012generalized}. Disjunctive constraint \eqref{eq:1ORC} becomes
\begin{subequations} \label{eq:CAV1BINARY}
 \begin{align}
       &x_c(t_f^*) - x_1(t_f^*) + (1-B_1)M\geq \varphi v_1(t_f^*)+\epsilon,\\
       &x_1(t_f^*) - x_c(t_f^*)+B_1M\geq \varphi v_c(t_f^*)+\epsilon.
\end{align}
\end{subequations}
Similarly, the other disjunctive constraints \eqref{eq:2ORC} and \eqref{eq:3ORC} can be converted. The extension to the arbitrary set of cooperating CAVs is as follows: 
\begin{subequations} \label{eq:ConstraintsFastLaneCAVC}
  \begin{align}
    &\forall j\in\mathcal{S}(t):\nonumber\\ 
       &x_C(t_f^*) - x_i(t_f^*) + (1-B_i)M\geq \varphi v_i(t_f^*)+\epsilon,\\
       &x_i(t_f^*) - x_C(t_f^*)+B_iM\geq \varphi v_c(t_f^*)+\epsilon,
\end{align}
\end{subequations}
\begin{comment}
    Similarly, disjunctive constraints \eqref{eq:2ORC} and \eqref{eq:3ORC} are converted as follows: 
\begin{subequations}  
  \begin{align}
       x_c(t_f^*) - x_2(t_f^*) + (1-B_2)M\geq& \varphi v_2(t_f^*)+\epsilon,\\
       x_2(t_f^*) - x_c(t_f^*)+B_2M\geq& \varphi v_c(t_f^*)+\epsilon.\\
       x_c(t_f^*) - x_3(t_f^*) + (1-B_3)M\geq& \varphi v_2(t_f^*)+\epsilon,\\
       x_2(t_f^*) - x_c(t_f^*)+B_3M\geq& \varphi v_c(t_f^*)+\epsilon.
\end{align}
\end{subequations}
\end{comment}

To further clarify the notation, $B_i = 1$ implies that CAV $C$ is in front of CAV $i\in \mathcal{S}(\tf)$. If $B_i =  1$, we need to have $B_{j}=1$ for $j\in\{i+1, \dots, m\}$, conveying that if CAV $C$ is in front of CAV $i$, it must be also in front of other CAVs located behind CAV $i$.
%%%%%%
\subsubsection{Optimization problem}
Next, we formally define the optimization problem for an arbitrary set of cooperating CAVs denoted by $\mathcal{S}(t_f^*)=\{1,\dots,m\}$ aimed at minimizing the global disruption of the maneuver on the fast lane.  

\begin{problem} \label{problem:P1}The following optimization problem is aimed at deriving the final states of cooperative CAVs on the fast lane, given the final time, $t_f^*$ and final states of CAV $C$ 
    \begin{subequations} 
    \allowdisplaybreaks
    \label{eq:optimalPair}
\begin{align}
        &\min_{\mathbf{x}_i(t_f^*),~B_i~\forall i\in \mathcal{S}(t_f^*)} \sum_{i\in \mathcal{S}(t_f^*)}  D_i(t_f^*) \nonumber\\
       % &\rm{s.t.}\quad \nonumber \\
        %\eqref{vconstraint}, \eqref{eq:ConstraintsFastLane}, \eqref{eq:reachable set}, \eqref{eq:fFrontc},  \eqref{eq:CFrontb}, \eqref{eq:ConstraintsFastLaneCAVC} \nonumber \\
         &x_F(t_f^*)-x_1(t_f^*)\geq\varphi v_1(t_f^*)+\epsilon,\\
         &x_{\color{black}{m}}(t_f^*)-x_B(t_f^*)\geq\varphi v_B(t_f^*)+\epsilon,\\
        &x_F(t_f^*) - x_C(t_f^*)\geq \varphi v_c(t_f^*)+\epsilon,\\
         &x_C(t_f^*) - x_B(t_f^*)\geq \varphi v_B(t_f^*)+\epsilon,\\
        &\forall j\in\mathcal{S}(t)\setminus\{m\}: \\
         &x_j(t_f^*)-x_{j+1}(t_f^*)\geq\varphi v_{j+1}(t_f^*)+\epsilon,\\
         &\forall i\in\mathcal{S}(t): \nonumber\\
       &x_C(t_f^*) - x_i(t_f^*) + (1-B_i)M\geq \varphi v_i(t_f^*)+\epsilon,\\
       &x_i(t_f^*) - x_C(t_f^*)+B_iM\geq \varphi v_c(t_f^*)+\epsilon,\\
       &p^{\rm upper}(x_i, v_i,t_f^*) \leq 0  \cup p^{\rm lower}(x_i, v_i,t_f^*) \geq 0,\\
       &0\leq v_{\min}\leq v_i(t_f^*)\leq v_{\max}.
    \end{align}
\end{subequations}
\end{problem}
\subsection{Unified Coordination Framework}\label{PA2}
In this section, we seek to eliminate the decoupling of time and space in the optimization setup.
Instead of deriving $t_f^*$ first and then solving for the final states of CAVs in the fast lane, this method seeks to solve for both simultaneously and to combine \ref{step1} and \ref{step2} (bottom green block in Fig. \ref{fig:blocks}).
\begin{problem}\label{problem:P2} The following optimization problem derives the final states of an arbitrary set of cooperative CAVs on the fast lane, final states of the ego CAV, and the final time, $t_f^*$.
     \begin{subequations}  
     \allowdisplaybreaks
     \label{eq:Unified}
\begin{align}
&\min_{t_f,\mathbf{x}_c(t_f),\mathbf{x}_i(t_f),B_i~\forall i\in \mathcal{S}(t_f)} \gamma_t t_f+D_c(t_f)+\sum_{i\in \mathcal{S}(t_f)}  D_i(t_f) \nonumber\\
        %&{\rm{s.t.}} \nonumber \\
        %& \eqref{vconstraint}, \eqref{eq:ConstraintsFastLane}, \eqref{eq:reachable set}, \eqref{eq:fFrontc},  \eqref{eq:CFrontb}, \eqref{eq:ConstraintsFastLaneCAVC} \nonumber \\
         &x_F(t_f)-x_1(t_f)\geq\varphi v_1(t_f)+\epsilon,\\
         &x_m(t_f)-x_B(t_f)\geq\varphi v_B(t_f)+\epsilon,\\
         &x_F(t_f) - x_c(t_f)\geq \varphi v_c(t_f)+\epsilon,\\
         &x_c(t_f) - x_B(t_f)\geq \varphi v_B(t_f)+\epsilon,\\
        &x_U(t_f) - x_C(t_f)\geq \varphi v_c(t_f)+\epsilon,\\
         &\forall j\in\mathcal{S}(t)\setminus\{m\}: \nonumber\\
         &x_j(t_f)-x_{j+1}(t_f)\geq\varphi v_{j+1}(t_f)+\epsilon,\\
         &\forall i\in\mathcal{S}(t): \nonumber\\
       &x_C(t_f) - x_i(t_f) + (1-B_i)M\geq \varphi v_i(t_f)+\epsilon,\\
       &x_i(t_f) - x_c(t_f)+B_iM\geq \varphi v_c(t_f)+\epsilon,\\
        &\forall i\in\mathcal{S}(t)\cup\{C\}: \nonumber\\
       &p^{\rm upper}(x_i, v_i,t_f) \leq 0  \cup p^{\rm lower}(x_i, v_i,t_f) \geq 0,\\
       &0\leq v_{\min}\leq v_i(t_f)\leq v_{\max},
    \end{align}
\end{subequations}   
where $\gamma_t$ is a weight to non-dimensionalize  the final time $t_f$.
\end{problem}
Note that Problems \ref{problem:P1} and \ref{problem:P2} are solved by CAV $C$. Next, the maneuver time and terminal states are broadcast to the CAVs on the fast lane to generate their optimal trajectories.  
\subsection{Trajectory Planner}\label{sec:trajPlanner}
The optimization problems detailed in equations \eqref{eq:optimalPair} and \eqref{eq:Unified} render the terminal states of all CAVs cooperating in the lane-change maneuver. 
However, to execute the maneuver, we also need to determine the trajectory plan over the specified time interval $[t_0, \tf)$, given the initial states $\mathbf{x}_i (t_0)$ and terminal states $\mathbf{x}_i (\tf)$ of vehicle $i$.
For clarity of notation in the formulation below, let the optimized terminal position and speed of vehicle $i$ be referred to as $x_i^{f}$ and $v_i^{f}$ respectively. 
% Note that $x_i^{\tf}$ and $v_i^{\tf}$ for each cooperating vehicle were found by solving the optimizations in  \eqref{eq:optimalPair} and \eqref{eq:Unified}.

\begin{subequations}
\allowdisplaybreaks
\label{eq: trajPlan}
    \begin{align}
        \min_{u_i(t)} & \int_{t_0}^{\tf} \cfrac{1}{2} u_i^2 (t) dt \\
        &{\textrm{s.t. \eqref{eq:dynamics}, \eqref{uconstraint}, \eqref{vconstraint}}}, \nonumber\\
        x_{i-1}(t) - x_i(t) \geq & \, \varphi v_i(t) + \epsilon, \, \forall t \in [t_0, \tf], \\
        (x_i(\tf) - x_i^{f})^2 \leq \, {\delta}_x, \quad&
        (v_i(\tf)  - v_i^{f})^2 \leq \, {\delta}_v,
    \end{align} 
\end{subequations}
where ${\delta}_x$ and ${\delta}_v$ are tolerances chosen for numerical feasibility.

The trajectories of all the vehicles can be determined sequentially in order in a distributed manner, i.e., $i = 1,2,\hdots$ and so on until the last vehicle in the cooperative set.
This allows passing the information of the preceding vehicle's trajectory to the current vehicle so as to ensure the satisfaction of safety constraints.

% \todoall{Todos for code}
% \textcolor{red}{\begin{enumerate}
%     \item \sout{Create a new branch. DONE}
%     \item \sout{Fix weights factors for disruption metric.}
%     \item \sout{Trajectory planner for PA \#1,2. Use sequential approach for now.}
%     \item \sout{IDM for other cars - only for previous phase papers. }
%     \item \sout{Write script to report terminal time, computation time, iterations required for feasibility, local and global disruptions, fuel consumption.}
%     \item \sout{ Vary number of cars, other initial conditions.}
% \end{enumerate}}

% {\color{blue} Not a to-do, but the previous approaches (conf, jnl.) don't impose jerk restrictions on the vehicles explicitly. It is problematic for us as well. 
% The ITSC (conf.) paper mentions it as future work, the journal paper doesn't mention it at all.}

\section{Simulation Results}\label{sec:simulation}

This section provides a summary of the simulation setup and the corresponding results demonstrating the various ways in which the proposed formulations outperform the baselines. 
The simulations were developed entirely in MATLAB with the help of CasADi \cite{andersson2019casadi} for numerical optimization on an Intel Core i7-1185G7 3.0 GHz.

The setting consists of a straight two-lane highway and an allowable speed range of $v = [5, 35]$ $m/s$. 
The headway parameter $\varphi$ and the safe distance parameter $\epsilon$ were chosen to be $0.2 \, s$ and $10 \, m$ respectively. 
% {\color{red} [$\epsilon$ may be too large. 1.5 m. in T-ITS paper.]}
We consider the case that all vehicles on the fast-lane are CAVs (i.e., vehicle $B$ and $F$ are not present). Congestion is generated due to the presence of a uncontrolled slow vehicle $U$ in the right lane traveling at a constant speed of $v_U = 20$ $m/s$.
CAV $C$ is initially present right behind $U$ and has a speed of $23 \, m/s$ as it seeks to change its lane. 
The control limits specified for every CAV are $u_{\rm min} = -7 \, m/s^2$ and $u_{\rm max} = 3.3 \, m/s^2$.
We used IPOPT \cite{wachter2006implementation} for obtaining the solutions for the optimal control problem in \eqref{eq: trajPlan}. The MINLP in \eqref{eq:optimalPair} and \eqref{eq:Unified} were solved using BONMIN \cite{bonami2011algorithms}. 

The desired flow speed was drawn from a uniform distribution to allow a range $v_{d} = [25,35] \, m/s$. The maximum time allowed for the maneuver was capped at $T_{\rm max} = 20 \, s$. Position and flow disruption were weighted equally, i.e., $\gamma = 0.5$, and $\gamma_t$ is chosen to be $\frac{1}{T_{\max}}$. 
For the trajectory planners, the allowable terminal state discrepancies were ${\delta}_x = 0.1 \, m^2$ and ${\delta}_v = 0.1 \, m^2/s^2$.
Once the optimal pair of cooperative vehicles were found for methods \cite{chavez2022sequential} and  \cite{armijos2022cooperative}, the vehicles in front followed a constant speed trajectory and the ones behind followed an IDM model \cite{treiber2000congested}.

\begin{table}
\renewcommand{\arraystretch}{1.15}
\caption{\label{results} Comparison of Averaged Results for Different Methods and Different Number of Cooperative CAVs}
\begin{center}
\begin{tabular}{|c|c| c c c c c|} 
\hline
 $|\mathcal{S}|$ &Methods & $\tf$ (s) &  $\Delta_{i, i+1}$ & $D_{S}$ & $n_{iter}$ & $t^{avg}_{iter} (s)$ \\ % [0.5ex] 
   % & & $(s)$ &   &  &  & $(s)$ \\ % [0.5ex] 
\hline
\multirow{4}{*}{6} & \cite{chavez2022sequential} & $2.81$	& $0.1171$	& $0.4426$	& $4.3$&$0.8$
\\  \cline{2-7}
& \cite{armijos2022cooperative} & $1.72$ & 	$0.0607$	& $0.3298$	& $6.5$& $4.49$
\\ \cline{2-7}
 &\ref{PA1} & $1.69$	& $0.0579$ & $0.0696$ & $6.3$& $0.33$
 \\ \cline{2-7}
 &\ref{PA2} & $2.42$ & $0.0007$ & $0.01$ &  $1$& $0.85$
 \\ \hline
 \hline
\multirow{4}{*}{5} & \cite{chavez2022sequential} & $2.67$&	$0.1219$	& $0.3949$ &	$3$&$0.77$
\\ 
 \cline{2-7}
& \cite{armijos2022cooperative} &$2.02$ &	$0.1015$ & 	$0.3394$ & 	$8.2$& $3.81$
\\
 \cline{2-7}
 &\ref{PA1} & $2.02$	&$ 0.1021$	& $0.1237$ &	$8.2$& $0.27$
 \\
 \cline{2-7}
 &\ref{PA2} & $2.56$	&$0.0004$ & $0.01$	& $1$&$0.76$\\
\hline \hline
\multirow{4}{*}{4} &  \cite{chavez2022sequential} &  $2.55$ & $0.0823$	& $0.1730$	& $3.2$& $1$
\\ 
 \cline{2-7}
& \cite{armijos2022cooperative} &$1.91$	&$ 0.0708$	& $0.1730$	& $7.3$&$3.17$
\\
 \cline{2-7}
 &\ref{PA1} & $1.7$	& $0.0329$	& $0.0389$	& $6.3$&$0.24$
 \\
 \cline{2-7}
 &\ref{PA2} & $2.25$ & $0.0004$ & $0.0076$	& $1$& $0.61$
 \\ 
 \hline \hline
 \multirow{4}{*}{3} & \cite{chavez2022sequential} & $2.57	$&$0.0692	$&$0.1145	$&$2.5$& $0.94$
\\ 
 \cline{2-7}
& \cite{armijos2022cooperative} &$1.92	$&$0.0568$	&$0.0894$	&$7.5$& $2.22$
\\
 \cline{2-7}
 &\ref{PA1} & $1.84	$&$0.0237$	&$0.0283$	&$7.2$& $0.18$
 \\
 \cline{2-7}
 &\ref{PA2} & $2.31	$&$0.0006$	&$0.0082$	&$1$& $0.5$
 \\ 
 \hline
\end{tabular}
\end{center}
\end{table}

The results of the simulations are summarized in Table \ref{results}. 
We compare our performance against the methods presented in \cite{chavez2022sequential} and \cite{armijos2022cooperative}.
We performed $24$ different simulations with random seeds varying initial conditions and free flow speeds to capture the impact of our proposed formulations on the disruption metrics under different number of cooperative CAVs varying from $3$ to $6$.
Here, $\tf$ is the terminal time found at the end of any time relaxation iterations carried out to ensure feasibility. 
$\Delta_{i, i+1}(\tf)$ refers to the local disruption incurred by the optimal pair of cooperating vehicles that CAV $C$ merges between, and $D_{\mathcal{S}} (\tf)$ denotes the global disruption that the full set of cooperating CAVs $\mathcal{S}(\tf)$ faces. 
% Both of these are non-dimensional quantities. 
$n_{Iters}$ refers to the average number of iterations that were required to ensure a feasible maneuver. %Lastly, $t_{iter}^{avg}$ refers to the mean computation time per iteration.

As can be clearly seen, the proposed approaches fare better than methods \cite{chavez2022sequential, armijos2022cooperative} in terms of both the local and global disruptions incurred due to the lane change maneuver.
The formulation described in section \ref{PA1} performs marginally better than the baselines and demonstrates the benefit of accounting for the entire set of cooperating vehicles in the fast lane. 
Moreover, the improvement in disruption by using the unified approach (\ref{PA2}) is orders of magnitude higher than any of the other methods. 
It also illustrates the redundancy of an iterative approach, ensuring feasibility in the very first optimization problem solved. To investigate the computational complexity of Problem \ref{problem:P1} and \ref{problem:P2}, we compute the mean of computation time per iteration (yellow blocks and green blocks in Fig. \ref{fig:blocks}) for each simulation of every approach, and averaged this means across all the seeds which are listed in Table \ref{results} as $t_{iter}^{avg}$. It shows that our approach is computationally feasible, however, the maximum number of CAVs to collaborate with should be upper-bounded based on the computational capabilities.

\section{Concluding Remarks and Discussion}\label{sec:Conclusion}

In this letter, we presented two formulations to perform a lane change maneuver for CAVs on a congested highway while minimizing a measure of disruption on the cooperating vehicles. 
The first approach seeks to account for the global impact on all the CAVs in the neighborhood instead of just the optimal pair. 
The second approach extends the first one by combining the multi-step optimization formulations to eliminate the decoupling between finding CAV $C$'s optimal trajectory and checking its feasibility for the cars on the fast lane, thereby removing the need for an iterative procedure.
Both methods perform better in terms of the disruption incurred by the cooperating CAVs for participating in the maneuver, and are a step towards designing more socially compliant formulations to relieve traffic congestion. 
Future studies will involve interactions with human-driven vehicles in a mixed autonomy framework. 

\bibliographystyle{IEEEtran.bst} 
\bibliography{reference/ref.bib}

\end{document}